\documentclass[11pt,a4paper]{article}
\usepackage[T1]{fontenc}
\usepackage{lmodern}
\usepackage[margin=26mm,headheight=15pt]{geometry}
\usepackage{amsmath,amssymb,amsthm,mathtools}
\usepackage{microtype,xcolor}
\usepackage[colorlinks=true,linkcolor=blue!45!black,citecolor=blue!45!black,urlcolor=blue!45!black]{hyperref}
\newtheorem{theorem}{Theorem}[section]
\newtheorem{lemma}[theorem]{Lemma}

\newtheorem{question}[theorem]{Question}

\usepackage{setspace}

\theoremstyle{definition}\newtheorem{definition}[theorem]{Definition}
\theoremstyle{remark}
\newcommand{\C}{\mathbb C}
\newcommand{\Cu}{C_u^*}
\newcommand{\Cum}{C_{u,\max}^*}
\newcommand{\norm}[1]{\lVert#1\rVert}

\DeclareMathOperator{\girth}{girth}
\DeclareMathOperator{\prop}{prop}
\title{Property A and the identification of the maximal and reduced uniform Roe algebras}

\author{Jintao Deng \and Geng Tian}
\date{}
\begin{document}
\maketitle
\onehalfspacing

\begin{abstract}
Let $X$ be a discrete metric space of bounded geometry. We prove that
$X$ has property A if and only if the canonical map 
$$\lambda:C^*_{u,max}(X)\to C^*_u(X)$$
is an isomorphism.
\end{abstract}

\section{Introduction}

Roe algebras arose in J.~Roe's study of higher indices of elliptic
operators on open manifolds \cite{Roe88}. The reduced uniform
Roe algebra appears in \cite[Section 4.4]{Roe03}, and its maximal
counterpart follows the construction of the maximal Roe algebra by
Gong, Wang and Yu \cite[Sections 3--4]{GWY}.

For a discrete metric space $X$ of bounded geometry, the identity on
the algebraic uniform Roe algebra extends to the canonical quotient map
$$
\lambda:\Cum(X)\longrightarrow\Cu(X).
$$
Whether $\lambda$ is an isomorphism is related to the large-scale
geometry of $X$. Yu introduced property A as a coarse analogue of
amenability \cite{Yu}. \v{S}pakula and Willett proved that property A
implies that $\lambda$ is an isomorphism and asked whether the converse
holds \cite[Proposition 1.3 and Section 1]{SW}. They also showed that
$\lambda$ induces an isomorphism on $K$-theory when $X$ coarsely embeds
into Hilbert space \cite[Theorem 1.4]{SW}. The box space of $\mathbb F_2$
constructed by Arzhantseva, Guentner and \v{S}pakula admits such an
embedding and fails property A \cite{AGS}. For that space, \v{S}pakula
and Willett explicitly asked whether $\lambda$ itself is an isomorphism
\cite[Example 1.16]{SW}.

\begin{question}\label{mainquestionofthispaper}
Let $X$ be a discrete metric space with bounded geometry. If the canonical
quotient map $\lambda$ is an isomorphism, must $X$ have Yu's property A?
\end{question}

Skandalis, Tu and Yu associated the coarse groupoid $G(X)$ to $X$
and recast the coarse assembly map as a groupoid Baum--Connes assembly
map with coefficients \cite[Proposition 4.8]{STY}. The maximal and
reduced uniform Roe algebras identify with the corresponding groupoid
$C^*$-algebras \cite[Remark 1.5]{SW}
$$
\Cum(X)\cong C^*_{\max}(G(X)),\qquad
\Cu(X)\cong C^*_r(G(X)).
$$
Under these identifications, $\lambda$ is the canonical groupoid
quotient. Skandalis, Tu and Yu also showed that $X$ has property A
if and only if $G(X)$ is topologically amenable
\cite[Theorem 5.3]{STY}. The question thus asks whether weak
containment forces topological amenability for coarse groupoids.

Willett's 2015 paper grew out of an attempt to answer this question
for the box space of $\mathrm{SL}(2,\mathbb Z)$ in \v{S}pakula and
Willett's Example~1.15 \cite[Remark 3.5]{Willett15}.
He constructed a Higson--Lafforgue--Skandalis groupoid with finite
quotients of $\mathbb F_2$ as its finite fibers and $\mathbb F_2$
as the fiber at infinity. The non-amenability of the fiber at
infinity forces the groupoid to fail topological amenability
\cite[Lemma 2.4]{Willett15}. Using Lubotzky and Shalom's property
$\mathrm{FD}$, Willett showed that the regular representations of
the chosen finite quotients recover the full norm on
$\mathbb C[\mathbb F_2]$. This yields equality of the maximal and
reduced groupoid $C^*$-algebras \cite[Lemmas 2.7 and 2.8]{Willett15}.
In Remark~3.5 of the same paper, he asked whether weak containment
can also hold for a non-amenable coarse groupoid.

Alekseev and Finn--Sell subsequently constructed non-amenable
principal \'etale groupoids with weak containment and second-countable
unit spaces \cite[Section 2]{AFS}. The coarse groupoid of an infinite
$X$ has the non-second-countable Stone--\v{C}ech compactification
$\beta X$ as its unit space \cite[Section 3]{STY}. Ara, Li, Lled\'o
and Wu recorded the coarse groupoid question as open
\cite[Example 4.6]{ALLW}. For coarse disjoint unions of large-girth
expanders, Buss, Echterhoff and Willett distinguished the two
uniform Roe completions by comparing the spectra of the graph
Laplacian \cite[proof of Lemma 4.7]{BEW}. 
We answer Question~\ref{mainquestionofthispaper} affirmatively.
Together with \v{S}pakula and Willett's result
\cite[Proposition 1.3]{SW}, this gives the following main theorem.

\begin{theorem}\label{thm:equivalence}
Let $X$ be a discrete metric space with bounded geometry. Then, $X$ has Yu's
property A if and only if the canonical quotient map
$$
\lambda:\Cum(X)\to\Cu(X)
$$
is an isomorphism.
\end{theorem}

The proof was inspired by the work of Roe--Willett \cite{RW}, in particular about the construction of the sequence of graphs using the failure off property A.  

When $X$ fails property A, we use Elek's result on uniform local
amenability \cite{Elek} to construct locally expanding graphs. We
then construct large-girth graphs on subsets of their vertex sets,
sign the edges using the theorem of Marcus, Spielman and Srivastava
\cite{MSS}, and form the signed Laplacian $L_\sigma$. This gives an
operator $T\in\Cum(X)$ with
$$
\norm{T}_{\max}=1,\qquad \norm{\lambda(T)}<1.
$$
Together with \v{S}pakula and Willett's forward implication
\cite[Proposition 1.3]{SW}, this yields our main result.

This paper is orgnized as follows. In Section~\ref{sec:roe}, we brielfy recall the two uniform Roe completions, Yu's
property A, and the decomposition of finite-propagation operators
into weighted partial translations. In Section~\ref{sec:graphs}, we
introduces the graph terminology and constructs large-girth graphs
from local expansion. Section~\ref{sec:state} constructs a state on
the maximal completion from signed paths. Section~\ref{sec:main-proof}
establishes the norm estimates and proves
Theorems \ref{thm:equivalence}.

\emph{\bf AI disclosure.}  We thank ChatGPT-6 for assistance in understanding the result of
Marcus, Spielman and Srivastava \cite{MSS} and for language editing
during the preparation of this manuscript.
All main results and mathematical proofs were developed by the authors.

\section{Uniform Roe algebras}\label{sec:roe}
In this section, we first briefly recall the definitions of maximal and reduced uniform Roe algebras and the concept of Yu's property A.

\subsection{Uniform Roe algebras}
Let $X$ be a discrete metric space with bounded geometry. We write
$B_X(x,R)=B_R(x)=\{z\in X:d(x,z)\le R\}$ and set
$$
N_X(R)=\sup_{x\in X}|B_X(x,R)|<\infty, \qquad R>0.
$$

Let $H_X=\ell^2(X)$. A bounded linear operator $T\in B(\ell^2(X))$ can be expressed as an $(X\times X)$-matrix
$$T=(T_{x,y})_{x,y\in X},$$
where $T_{x,y}\in \mathbb{C}$ for any $x,y\in X$. The \emph{propagation} of $T$ is defined as 
$$
\prop(T)=\sup\{d(x,y):T_{x,y}\neq 0\},\qquad \prop(0)=0.
$$

We say that the operator $T$ has \emph{finite propagation} if 
$$\prop(T)<\infty.$$

Now let us recall the definitions of the algebraic uniform Roe algebras.
\begin{definition}
Let $X$ be a discrete metric space with bounded geometry. The algebraic uniform Roe algebra, denoted by $\mathbb{C}_u[X]$, is the $*$-algebra of all bounded operators $T\in B(\ell^2(X))$ with finite propagation.    
\end{definition}

For any two operators $T=(T_{x,y})_{x,y\in X}, S=(S_{x,y})_{x,y\in X}\in \mathbb{C}_u[X]$, the product $TS$ is defined as 
$$(TS)_{x,z}=\sum_{y\in X}T_{x,y}S_{y,z},$$
for any $x,z \in X$.

For every $r>0$, there exists a constant $C_r>0$ such that 
$$\|\pi(T)\|\leq C_r \cdot \sup_{x,y\in X} |T_{x,y}|$$
for any $*$-representation $\pi: \mathbb{C}_u[X]\to B(H)$ and
for any $T\in \mathbb{C}_u[X]$ with $\prop(T)\leq r$ \cite[Lemma 1.10]{SW}. As a result, the following maximal norm is finite: 
$$
\|T\|_{\max}=\sup \{\|\pi(T)\|: \pi: \mathbb{C}_u[X]\to B(H)~\mbox{is a $*$-}\mbox{representation}\}.
$$
The maximal and reduced uniform Roe algebras can be viewed as coarse
analogues of the maximal and reduced group $C^*$-algebras.
\begin{definition}[{\cite[Definitions 1.9 and 1.11]{SW}}]
Let $X$ be a discrete metric space with bounded geometry.
    \begin{enumerate}
        \item[(1)] The uniform Roe algebra $C_u^*(X)$ is defined to be the completion of $\mathbb{C}_u[X]$ under the norm on $B(\ell^2(X))$;
        \item[(2)] The maximal uniform Roe algebra $C^*_{u,max}(X)$ is defined to be the completion of $\mathbb{C}_u[X]$ under the maximal norm.
    \end{enumerate}
\end{definition}

For the index-theoretic origins of Roe algebras, see \cite{Roe88}.
In Roe's terminology, $C_u^*(X)$ is the
\emph{uniform translation $C^*$-algebra};
see \cite[Section 4.4]{Roe03}.

By the universality of the maximal norm, the identity map on $\mathbb{C}_u[X]$ continuously extends to the quotient map 
$$
\lambda: C^*_{u,max}(X)\to C^*_{u}(X).
$$
We use $\norm{\cdot}_{\mathrm{red}}$ for the concrete operator norm.

We use the following equivalent formulation of Yu's property A \cite{Yu}.
\begin{definition}[{\cite[Definition 2.1]{SW}}]\label{def:property-a}
    Let $X$ be a discrete metric space with bounded geometry, and let $\ell^2(X)_+$ denote the nonnegative real-valued functions in $\ell^2(X)$. We say that $X$ has Yu's property A if for every $R>0$ and $\epsilon>0$, there exists a map $\xi: X\to \ell^2(X)_+$ satisfying
    \begin{enumerate}
        \item[(1)] for each $x\in X$, $\|\xi_x\|_2=1$;
        \item[(2)] there exists an $S>0$, such that ${\rm Supp}(\xi_x)\subseteq B_S(x)$, for all $x\in X$;
        \item[(3)] for all $x, y\in X$ with $d(x,y)\leq R$, $\|\xi_x-\xi_y\|_2<\epsilon$.
    \end{enumerate}
\end{definition}

The following result of \v{S}pakula and Willett gives the forward implication of Theorem~\ref{thm:equivalence}.
\begin{theorem}[{\cite[Proposition 1.3]{SW}}]
Let $X$ be a discrete metric space with bounded geometry. If $X$ has Yu's property A, then the quotient map 
    $$
\lambda: C^*_{u,max}(X)\to C^*_{u}(X)
$$
is an isomorphism.
\end{theorem}

\subsection{Partial translations}\label{sec:partial-translations}

A \emph{partial translation} is a bijection $t:A\to B$ between subsets
of $X$ with $\sup_{x\in A}d_X(x,t(x))<\infty$. Its operator is
$$
V_t\delta_x=\begin{cases}\delta_{t(x)},&x\in A,\\0,&x\notin A.\end{cases}
$$
For $f\in\ell^\infty(X)$, write $M_f$ for diagonal multiplication by $f$.

Every finite-propagation operator is a finite sum of operators
$M_f V_t$. Indeed, if the propagation is at most $R$, each row and
column contains at most $N_X(R)$ nonzero entries. Enumerate the
nonzero entries and color them greedily with $2N_X(R)-1$ colors,
forbidding repetitions of a color in the same row or column.
At most $2(N_X(R)-1)$ colors are forbidden at each step. Within each
color, the nonzero entries define a partial translation from column
indices to row indices, and their values define the coefficient $f$
on its range.

\section{Local expansion and large-girth graphs}\label{sec:graphs}

Let $X$ be a discrete metric space without property A. We first construct a squence of finite graphs with uniform local expansion. Based on this sequence, we construct a sequence of graphs
with degrees five or six, uniformly bounded ambient edge lengths, and
girth tending to infinity.

\subsection{Preliminaries on graphs}\label{sec:graph-preliminaries}

In this subsection, we recall some background on graph theory. 

A \emph{simple graph} is a pair $G=(V(G),E(G))$, where $V(G)$ is
its vertex set and $E(G)$ is a set of unordered pairs of distinct
vertices. We write $(x,y)$ for the edge with endpoints $x$ and $y$.
The expressions $(x,y)$ and $(y,x)$ denote the same undirected edge.
The endpoints of an edge are \emph{adjacent}, or \emph{neighbors}.
The neighbor set and degree of $x$ are
$$
N_G(x)=\{y:(x,y)\in E(G)\},\qquad d_G(x)=|N_G(x)|.
$$
An edge is \emph{incident} to each of its endpoints. A graph has
\emph{bounded degree} if there is a finite constant $\delta$ such that
$d_G(x)\le \delta$ for every vertex $x\in V(G)$.

A \emph{subgraph} $F$ of $G$ is a graph whose vertex set $V(F)$ and edge set $E(F)$ are subsets of $V(G)$ and $E(G)$, respectively, and every edge of $E(F)$ has endpoints in $V(F)$. For
$W\subseteq V(G)$, the \emph{induced subgraph} $G[W]$ has vertex set
$W$ and all edges of $G$ with both endpoints in $W$. A simple graph
\emph{on} $W$ has vertex set $W$, with an edge set specified separately.

A \emph{matching} is a set of edges no two of which share an
endpoint. A \emph{proper edge coloring} assigns colors to edges so
that edges sharing an endpoint have different colors. Each color
class is therefore a matching.

A \emph{walk of length $k$} is a finite sequence
$$
x_0,e_1,x_1,e_2,\ldots,e_k,x_k,
$$
where $e_j$ has endpoints $x_{j-1}$ and $x_j$. A walk is \emph{closed}
if $x_0=x_k$. A \emph{trail} is a walk with no repeated edge, and a
\emph{path} is a walk with pairwise distinct vertices. A graph is
\emph{connected} if any two vertices can be joined by a path. Its
\emph{components} are the maximal connected subgraphs. A \emph{cycle}
is a closed trail of positive length in which $x_0,\ldots,x_{k-1}$
are pairwise distinct. Cycles are identified under change of starting
vertex and reversal. Every closed trail of positive length contains
a cycle, by taking a shortest nonempty closed subtrail.

The \emph{girth} $\girth(G)$ is the least length of a cycle, with value
$\infty$ if there is no cycle. A \emph{forest} is a graph with no cycles.
A \emph{tree} is a connected forest. In a tree, there is a unique path
between any two vertices, since distinct paths with the same endpoints
would contain a cycle in their union.

For vertices $x,y$ in the same component, the \emph{graph distance}
$d_G^{\mathrm{gr}}(x,y)$ is defined to be the minimum length of a path joining them. For vertices in different components, we set $d_G^{\mathrm{gr}}(x,y)=\infty$.

If $V(G)\subseteq X$, the \emph{ambient length} of the edge
$(x,y)$ is $d_X(x,y)$. A bound $R$ on all ambient edge lengths
implies
$$
d_X(x,y)\le R\,d_G^{\mathrm{gr}}(x,y)
\quad\text{when $x,y$ are in the same component}.
$$
For nonempty $E\subseteq X$, we write
$\operatorname{diam}_X(E)=\sup\{d_X(x,y):x,y\in E\}$.

\begin{definition}\label{def:scale-graph}
Let $(X,d_X)$ be a discrete metric space with bounded geometry and $R>0$.
The \emph{scale-$R$ graph} $\Gamma_R$ has vertex set $X$ and an
edge $(x,y)$ for each pair of distinct vertices with $d_X(x,y)\le R$.
\end{definition}

\begin{definition}\label{def:graphboundary}
Let $G$ be a graph and $E\subseteq V(G)$. The \emph{internal vertex boundary} and \emph{edge boundary} are defined as
\begin{equation*}
\begin{split}
\partial_G^{\mathrm{int}}E
&=\{x\in E:\text{some }y\in V(G)\setminus E
                         \text{ satisfies }(x,y)\in E(G)\},\\
\partial_G^{\mathrm e}E
&=\{e\in E(G):e\text{ has exactly one endpoint in }E\}.
\end{split}
\end{equation*}
\end{definition}
Every internal boundary vertex is incident to a crossing edge, so $|\partial_G^{\mathrm{int}}E|\le|\partial_G^{\mathrm e}E|$.
If $F$ is a subgraph, $\partial_F^{\mathrm{int}}E$ counts only
neighbors in $V(F)\setminus E$ and edges belonging to $F$.

\subsection{Local expansion}

Brodzki, Niblo, \v{S}pakula, Willett and Wright introduced uniform
local amenability (ULA) and proved its coarse invariance
\cite[Definition 2.2 and Lemma 2.6]{BNSWW}. We use Elek's formulation
\cite{Elek,ElekTimar} for graphs with bounded degree.

\begin{definition}[{\cite[Definition 1.2]{Elek}}]\label{def:ULA}
Let $G$ be an infinite graph of bounded degree. The graph $G$ is said to be \emph{uniformly locally amenable} (ULA)
if, for every $h>0$, there is $K<\infty$ such that each nonempty finite subgraph $F$ contains a nonempty $E\subseteq V(F)$ with
$$
|E|\le K,\qquad |\partial_F^{\mathrm{int}}E|\le h|E|.
$$
\end{definition}

For metric spaces with bounded geometry, Brodzki, Niblo, \v{S}pakula, Willett and Wright showed that Yu's property A implies ULA
\cite[Proposition 3.2 and Lemma 2.7]{BNSWW}, while Elek proved the
converse \cite[Theorem 1]{Elek}. For the equivalence in the
graph-metric formulation, see also \cite[Theorem 1]{ElekTimar}.

The following construction is inspired by the arguments in Roe--Willett \cite{RW}. Failure of ULA gives rise to a sequence of finite graphs with a uniform lower bound on the ratio of internal boundary size to cardinality for small vertex sets. 
\begin{lemma}\label{lem:witnesses}
Suppose $X$ has bounded geometry without Yu's property A. Then, there are
constants $R>0$, $h\in(0,1)$ and an integer $\delta\ge2$, and a sequence
$(G_n)_n$ of finite connected simple graphs with nonempty,
pairwise disjoint vertex sets $(W_n)_n$ in $X$, such that 

\begin{enumerate}
\item[(1)] each $G_n$ has degrees
at most $\delta$ and ambient edge lengths at most $R$;
\item[(2)] for every $n$ and every nonempty susbet
$E\subseteq W_n$ with
$\operatorname{diam}_X(E)\le2nR$, we have that
\begin{equation}\label{eq:boundary}
|\partial_{G_n}^{\mathrm{int}}E|>h|E|.
\end{equation}

\end{enumerate}
\end{lemma}

\begin{proof}
It is obvious that the space $X$ is infinite.
Suppose that every scale graph $\Gamma_R$ has ULA. Fix
$R,\varepsilon>0$. Since $\Gamma_R$ has bounded degree, ULA gives
nonnegative $\ell^1$ unit vectors $a_x$ supported in graph balls
of a common radius $S$, with
$\lVert a_x-a_y\rVert_1<\varepsilon^2$ on edges. Each such ball
about $x$ is contained in $B_X(x,RS)$, and distinct points at
ambient distance at most $R$ are adjacent in $\Gamma_R$. Thus
$\xi_x=\sqrt{a_x}$ gives the required $\ell^2$ unit vectors at
scale $R$, since
$\lVert\xi_x-\xi_y\rVert_2^2\le\lVert a_x-a_y\rVert_1$.
As $R$ and $\varepsilon$ are arbitrary, $X$ has property A,
a contradiction. Therefore, there must be an $R>0$, such that the scaled graph $\Gamma_R$ does not have ULA. 

Fix $R>0$ for which $\Gamma_R$ does not have ULA. There is $h_0\in(0,1)$
such that for every positive integer $K$ there is a nonempty finite
subgraph $F$ of $\Gamma_R$ satisfying
\begin{equation}\label{eq:smallsets}
|\partial_F^{\mathrm{int}}A|>h_0|A|
\qquad (A\subseteq V(F),\ 1\le |A|\le K).
\end{equation}
In particular, $|V(F)|>K$, since $V(F)$ has empty internal
boundary in $F$. Let 
$$h=h_0/2, ~\mbox{and}~\delta=\max\{2,N_X(R)-1\},$$
where $N_X(R)=\sup_{x\in X}|B_X(x,R)|$.

We construct $(G_n)_n$ by induction. We start with $G_1$.
Choose $F_1$ satisfying \eqref{eq:smallsets} with $K=N_X(2R)$
and take a connected component $G_1$ with vertex set $W_1$.
A nonempty subset $E\subseteq W_1$ of ambient diameter at most
$2R$ lies in $B_X(x,2R)$ for any $x\in E$, so $|E|\le K$.
Its internal boundary is the same in $G_1$ and $F_1$, so
\eqref{eq:smallsets} gives the required bound with constant $h$.

For the construction of $G_2$, take $M=W_1$ and $k=N_X(4R)$.
Choose an integer $b>|M|/h$ and a new $F_2$ satisfying
\eqref{eq:smallsets} with $K=b+k$. First remove
$S_0=V(F_2)\cap M$. Let $S$ be the set of vertices removed so far
and write $F_S=F_2[V(F_2)\setminus S]$ for the remaining induced
subgraph. Whenever $F_S$ contains a nonempty subset $E$ with
$$
|E|\le k,\qquad |\partial_{F_S}^{\mathrm{int}}E|\le h|E|,
$$
add $E$ to $S$ and repeat.

To see why this leaves a nonempty graph, count the internal
boundary of $S$ in the original graph $F_2$. When $E$ is added,
a vertex in $S$ can remain on the boundary only if it was already
there. The vertices of $E$ on the new boundary are exactly those
with a neighbor outside $S\cup E$, namely the vertices in
$\partial_{F_S}^{\mathrm{int}}E$. Thus each deletion adds at most
$h|E|$ boundary vertices. The initial boundary has size at most
$|M|$, and the sets removed at successive steps are disjoint.
It follows that
$$
|\partial_{F_2}^{\mathrm{int}}S|\le |M|+h|S|.
$$

Initially $|S|\le |M|<b$. If $|S|$ first reaches or exceeds $b$,
then $b\le |S|\le b+k=K$, since the last step adds at most $k$
vertices. Thus \eqref{eq:smallsets} applies to $S$. But
$|M|<hb\le h|S|$, so
$$
|\partial_{F_2}^{\mathrm{int}}S|
\le |M|+h|S|<2h|S|=h_0|S|,
$$
contrary to \eqref{eq:smallsets}. Hence $|S|<b$ throughout.
The process terminates since $F_2$ is finite and each step
removes at least one vertex. Since $|V(F_2)|>K$, the final
$F_S$ is nonempty.

Take a connected component $G_2$ of the final $F_S$ and put
$W_2=V(G_2)$. The initial deletion gives
$W_2\cap W_1=\varnothing$. If $\varnothing\ne E\subseteq W_2$
has ambient diameter at most $4R$, then $E\subseteq B_X(x,4R)$
for any $x\in E$, so $|E|\le k$. The stopping condition gives
$|\partial_{F_S}^{\mathrm{int}}E|>h|E|$. Every neighbor in $F_S$
of a vertex of $E$ lies in $G_2$, so
$\partial_{G_2}^{\mathrm{int}}E=\partial_{F_S}^{\mathrm{int}}E$.

Now suppose that $G_1,\ldots,G_n$ have been constructed. Repeat
the construction of $G_2$ with $M=\bigcup_{i=1}^{n}W_i$ and
$k=N_X(2(n+1)R)$. The resulting component $G_{n+1}$ avoids $M$
and satisfies the required boundary estimate with the same
constant $h$. Each $G_n$ is a subgraph of $\Gamma_R$, so the
degree bound $D$ and the ambient edge length bound $R$ hold
uniformly in $n$.
\end{proof}

\begin{definition}\label{def:graphoperators}
For a finite simple graph $G$, the \emph{adjacency matrix} $A_G$
acts on $\ell^2(V(G))$ and is defined by
$$
(A_G)_{xy}=\begin{cases}1,&(x,y)\in E(G),\\0,&\text{otherwise}.
\end{cases}
\qquad (A_G f)(x)=\sum_{y\in N_G(x)}f(y).
$$
The \emph{degree matrix} is
$D_G=\operatorname{diag}(d_G(x))$, and the \emph{combinatorial
Laplacian} is $\Delta_G=D_G-A_G$. 
\end{definition}
With each undirected edge
counted once, its quadratic form is
$$
\langle f,\Delta_G f\rangle
=\sum_{(x,y)\in E(G)}|f(x)-f(y)|^2.
$$
Let $(G_n)_n$ be the graphs with vertex sets $(W_n)_n$ in Lemma~\ref{lem:witnesses}. On $\ell^2(W_n)$,  consider the following operators
$$
\Delta_n=D_{G_n}-A_{G_n},\qquad
P_n=I-\frac{\Delta_n}{2\delta}.
$$
Each $P_n$ is symmetric and positive semidefinite, with nonnegative
entries and row sums one.
Consider the constant
$$\rho=1-\frac{h^2}{4\delta^2}\in(0,1).$$

For the operator $P_n$, we have the following estimates on the diagonals of its powers $P_n^{\ell}$, for $\ell\leq n$.

\begin{lemma}\label{lem:local-matrix-bound}
For $x\in W_n$ and $1\le \ell\le n$,
\begin{equation}\label{eq:local-matrix-bound}
(P_n^\ell)_{xx}\le\rho^\ell.
\end{equation}
\end{lemma}
\begin{proof}

Fix $x\in W_n$ and $1\le\ell\le n$. Let
$B=W_n\cap B_X(x,\ell R)$ and let $P_{n,B}$ be the compression
of $P_n$ to $\ell^2(B)$. For any function $f\in\ell^2(B)$, it extends by zero to $W_n$.
List the distinct values of $|f|^2$ as
$t_1<\cdots<t_s$, put $t_0=0$, and set
$$E_i=\{u:|f(u)|^2\ge t_i\}.$$ 
Since $\operatorname{diam}_X(B)\le 2nR$, by the Lemma~\ref{lem:witnesses}, we have that $|\partial_{G_n}^{\mathrm e}E_i|\ge h|E_i|$. 
Summing
these inequalities with weights $t_i-t_{i-1}$ gives rise to the inequality
$$
h\|f\|_2^2
\le \sum_{i=1}^s(t_i-t_{i-1})|\partial_{G_n}^{\mathrm e}E_i|
= \sum_{(u,v)\in E(G_n)}
  \bigl||f(u)|^2-|f(v)|^2\bigr|.
$$
By Cauchy--Schwarz and the degree bound, the last sum is at most
$$
\left(\sum_{(u,v)\in E(G_n)}|f(u)-f(v)|^2\right)^{1/2}
\left(\sum_{(u,v)\in E(G_n)}(|f(u)|+|f(v)|)^2\right)^{1/2}
\le \langle f,\Delta_n f\rangle^{1/2}
    (2\delta\|f\|_2^2)^{1/2}.
$$
It follows that $\langle f,\Delta_n f\rangle
\ge \frac{h^2}{2\delta}\|f\|_2^2$, and hence
$$
0\le \langle f,P_{n,B}f\rangle
=\|f\|_2^2-\frac{1}{2\delta}\langle f,\Delta_n f\rangle
\le \rho\|f\|_2^2.
$$
Thus $\|P_{n,B}\|\le\rho$.

Each multiplication by $P_n$ enlarges the support by at most $R$
in the ambient metric. Hence $P_n^j\delta_x$ is supported in $B$
for $0\le j\le\ell$. Projecting onto $\ell^2(B)$ after each step
therefore leaves these vectors unchanged, so
$P_{n,B}^{\ell}\delta_x=P_n^{\ell}\delta_x$. Accordngly, we have that
$$
(P_n^\ell)_{xx}
=\langle\delta_x,P_{n,B}^{\ell}\delta_x\rangle
\le \|P_{n,B}\|^\ell
\le \rho^\ell.
$$
\end{proof}

We use the estimates on the entries of $P_n^m$ to define a sequence of new graph and then removing certain vertice and edges gives us a sequence of graphs with large grith and degree at most $6$.
\begin{lemma}\label{lem:extraction}
Let $(G_n)_n$ be as in Lemma~\ref{lem:witnesses}, with vertex sets
$(W_n)_n$. There exist a constant $L<\infty$, an integer $n_0$ and a sequence
$(H_n)_{n\ge n_0}$ of simple graphs with nonempty vertex sets
$(U_n)_{n\ge n_0}$ such that $U_n\subseteq W_n$ for each $n\ge n_0$ and
$$
d_{H_n}(x)\in\{5,6\}\quad(x\in U_n),\qquad
\girth(H_n)\longrightarrow\infty,\qquad
d_X(x,y)\le L\quad((x,y)\in E(H_n)).
$$
\end{lemma}
\begin{proof}
Fix an integer $m\ge 1$ with $\rho^m\le 1/20$. For $n\ge 4m$, denote
$$
Q=P_n^m,\qquad q_n=\lfloor n/m\rfloor+1,\qquad N=|W_n|.
$$
The matrix $Q$ is symmetric and positive semidefinite, with
nonnegative entries and row sums one. By Lemma~\ref{lem:local-matrix-bound},
$$
0\le Q_{xy}\le \sqrt{Q_{xx}Q_{yy}}\le \frac1{20}.
$$
Let $J$ be the simple graph on $W_n$ with edge set
$$
E(J)=\{(x,y):x\ne y,\ Q_{xy}>0\}.
$$
Since $P_n$ has propagation at most $R$, every edge of $J$ has
ambient length at most $mR$.

For a graph $K$, write $\mathcal C(K)$ for its cycles of length
less than $q_n$. For $u\in[0,1]^{E(J)}$, define
$$
F(u)=\sum_{e\in E(J)}u_e
-\sum_{C\in\mathcal C(J)}\prod_{e\in E(C)}u_e.
$$
Set $a_{(x,y)}=10Q_{xy}$, so that $a_e\in[0,1/2]$. By symmetry
and the row sums of $Q$,
$$
\sum_{e\in E(J)}a_e
=5\sum_{x\in W_n}(1-Q_{xx})\ge \frac{19}{4}N.
$$
Choose a starting vertex and a direction for each $k$-cycle of $J$.
Distinct cycles then give distinct terms in the expansion of
$\operatorname{tr}(Q^k)$. Each such term is the product of the
corresponding edge weights $a_e$ divided by $10^k$. Since the trace
expansion has nonnegative terms and $mk\le n$ for $k<q_n$, we obtain that
$$
\sum_{\substack{C\in\mathcal C(J)\\ |E(C)|=k}}
  \prod_{e\in E(C)}a_e
\le 10^k\operatorname{tr}(Q^k)
\le (10\rho^m)^kN
\le 2^{-k}N
\qquad (3\le k<q_n).
$$
Consequently,
$$
F(a)\ge
\left(\frac{19}{4}-\sum_{k=3}^\infty 2^{-k}\right)N
=\frac92N>4N.
$$

With the other coordinates fixed, $F$ is affine in $u_e$. Thus
one of the choices $u_e=0$ and $u_e=1$ preserves or increases $F$.
Applying this successively to all coordinates of $a$ gives
$b\in\{0,1\}^{E(J)}$ with $F(b)\ge F(a)$. Let $J_0$ be the graph
on $W_n$ with edges $\{e\in E(J):b_e=1\}$. Then
$$
|E(J_0)|-|\mathcal C(J_0)|=F(b)>4N.
$$
For each $C\in\mathcal C(J_0)$ choose $e_C\in E(C)$, and let $J_1$
be obtained by deleting the set $\{e_C:C\in\mathcal C(J_0)\}$.
Hence
$$
|E(J_1)|\ge |E(J_0)|-|\mathcal C(J_0)|>4N,
\qquad \operatorname{girth}(J_1)\ge q_n.
$$

Repeatedly delete vertices of current degree at most four from
$J_1$. If all vertices were deleted in the order $v_1,\ldots,v_N$,
counting each edge when its first endpoint is deleted would result in
$$
|E(J_1)|
=\sum_{i=1}^{N}d_{J_1[\{v_i,\ldots,v_N\}]}(v_i)
\le 4N,
$$
and this is a contradiction. The remaining graph $K$ therefore has a nonempty
vertex set $U_n\subseteq W_n$ and satisfies
$$
\min_{v\in U_n}d_K(v)\ge 5,\qquad
\operatorname{girth}(K)\ge q_n.
$$

Starting from $K$, delete edges whose two endpoints both have
current degree greater than five until no such edge remains.
Write $K'$ for the resulting graph and put
$$
S=\{v\in U_n:d_{K'}(v)>5\}.
$$
Then $S$ is independent and $d_{K'}(u)=5$ for $u\in U_n\setminus S$.
For each $v\in S$, partition $2r_v$ of its incident edges into pairs,
where
$$
r_v=\left\lfloor\frac{d_{K'}(v)-5}{2}\right\rfloor.
$$
Each pair forms a length-two path $u-v-w$ with distinct endpoints
in $U_n\setminus S$. The chosen paths are pairwise edge-disjoint,
since the pairs at each $v$ are disjoint and each edge of $K'$
meets at most one vertex of $S$.

Let $\mathcal P$ consist of the chosen length-two paths together
with the length-one paths given by the unused edges of $K'$.
For a path $p\in\mathcal P$ with endpoints $x,y$, write $\ell(p)$
for its length and set $\partial p=(x,y)$. Define
$$
V(H_n)=U_n,\qquad E(H_n)=\{\partial p:p\in\mathcal P\}.
$$
If distinct $p,p'\in\mathcal P$ had the same endpoints, traversing
$p$ and returning along $p'$ would give a closed trail in $K'$.
Since a nonempty closed trail contains a cycle, this would give
$$
q_n\le \operatorname{girth}(K')
\le \ell(p)+\ell(p')\le 4,
$$
contrary to $q_n\ge 5$. Thus the endpoint map is injective and
$H_n$ is simple.

For $u\in U_n\setminus S$, let $s_u$ be the number of chosen
length-two paths with endpoint $u$. Each such path replaces one
edge incident to $u$ by one shortcut. For $v\in S$, each of the
$r_v$ chosen paths removes two edges incident to $v$. Hence
$$
\begin{aligned}
d_{H_n}(u)&=d_{K'}(u)-s_u+s_u=5
&& (u\in U_n\setminus S),\\
d_{H_n}(v)&=d_{K'}(v)-2r_v\in\{5,6\}
&& (v\in S).
\end{aligned}
$$

Let $p_e$ be the unique path representing $e\in E(H_n)$.
For a cycle $C$ in $H_n$, concatenating the paths $p_e$ in cycle
order gives a closed trail in $K'$, since the paths are pairwise
edge-disjoint. Therefore
$$
q_n\le \sum_{e\in E(C)}\ell(p_e)\le 2|E(C)|,
\qquad
\operatorname{girth}(H_n)\ge \frac{q_n}{2}\longrightarrow\infty.
$$
Finally, if $e=(x,y)\in E(H_n)$, the triangle inequality along
$p_e$ gives
$$
d_X(x,y)\le mR\,\ell(p_e)\le 2mR.
$$
The conclusion follows with $n_0=4m$ and $L=2mR$.
\end{proof}
We would like to remark here that the shortcut construction follows from the arguments in \cite[Lemma~14]{Kun}.
\section{A state on the maximal uniform Roe algebra}\label{sec:state}

In this section, we construct a state on the maximal
uniform Roe algebra using products of edge signs along short paths, for any given signing of graphs.
The state vanishes on the signed Laplacian and provides the lower bound
for the maximal norm in Section~\ref{sec:main-proof}.

Let $(H_n)_n$ satisfy the conclusion of Lemma~\ref{lem:extraction},
with pairwise disjoint vertex sets $(U_n)_n$. The degrees are five or six,
the ambient edge lengths are at most $L$, and $\girth(H_n)\to\infty$.
After discarding a finite initial segment, we have $\girth(H_n)>6$
for every remaining $n$. Set $U=\bigsqcup_n U_n$ and $p=1_U$ for
this sequence.

A \emph{signing} of a graph $H$ is a map 
$$\sigma:E(H)\to\{-1,1\}.$$
Its \emph{signed adjacency matrix} has entry $\sigma((x,y))$ when
$(x,y)$ is an edge and zero otherwise. For a sequence of signings
$\sigma=(\sigma_n)_n$, let $D$ and $B_\sigma$ be the block
diagonal operators with blocks $D_{H_n}$ and the signed adjacency
matrices, respectively, extended by zero off $U$. The degree diagonal
has norm at most six. For $B_\sigma$, Cauchy--Schwarz bounds
$|(B_\sigma f)(x)|^2$ by six times the sum of $|f(y)|^2$ over
neighbors $y$ of $x$. Summing over $x$ counts each $|f(y)|^2$ at
most six times, so $\norm{B_\sigma}\le6$. The ambient edge length
bound gives propagation at most $L$. Define
$$
L_\sigma=D-B_\sigma,\qquad T=p-L_\sigma/12.
$$
The \emph{sign of a path} is the product of its edge signs, with
value one on the empty path. 

Since $\girth(H_n)\to\infty$, choose integers $J_n\ge1$
with $J_n\to\infty$ and $3\cdot 2^{J_n}<\girth(H_n)$.
Write $d_n=d_{H_n}^{\mathrm{gr}}$. For $1\le j\le J_n$, let $r_j=2^j$, and define
$$
K_{n,j}(x,y)=
\begin{cases}
\text{sign of the shortest path from $x$ to $y$}, & d_n(x,y)\le r_j,\\
0, & d_n(x,y)>r_j.
\end{cases}
$$
Two distinct paths of length at most $r_j$ with the same endpoints
would contain a cycle of length at most $2r_j$ in their union.
The girth bound therefore makes the sign in the first case
well defined. Define
$$
K_n=\frac{1}{J_n}\sum_{j=1}^{J_n}K_{n,j}.
$$
Then $K_n$ is real and symmetric, its entries have absolute value
at most one, $K_n(x,x)=1$, and $K_n(x,y)=\sigma_n((x,y))$ on
each edge of $H_n$.

\begin{lemma}\label{lem:small-kernel}
Let $M$ be a positive integer. For any $S\subseteq U_n$ with
$|S|\le M$, there is a positive semidefinite matrix on $S$ that
differs from $K_n|_{S\times S}$ entrywise by at most
$\binom{M}{2}/J_n$.
\end{lemma}
\begin{proof}
Let $I_S$ be the set of indices $j$ in the average defining $K_n$
for which no pair $x,y\in S$ satisfies
$r_j<d_n(x,y)\le 2r_j$. Since $2r_j=r_{j+1}$, the intervals
$(r_j,2r_j]$ are disjoint. Each unordered pair in $S$ therefore
excludes at most one index, so at most $\binom{M}{2}$ indices
lie outside $I_S$.

Fix $j\in I_S$. If $d_n(x,y)\le r_j$ and $d_n(y,z)\le r_j$,
then $d_n(x,z)\le 2r_j$ by the triangle inequality. By the definition
of $I_S$, we have that $d_n(x,z)\le r_j$. Thus the relation
$d_n(x,y)\le r_j$ on $S$ is an equivalence relation. Therefore,
$K_{n,j}(x,y)=0$ for vertices in distinct classes.

Fix a root $o$ in an equivalence class $C$ and write $s_x$ for the
sign of the shortest path from $o$ to $x\in C$. For $x,y\in C$,
the three shortest paths between $o,x,y$ have a connected union
with at most $3r_j<\girth(H_n)$ edges, so their union is a tree.
The signs on the common initial segment of the paths from $o$
to $x$ and from $o$ to $y$ cancel in $s_xs_y$, leaving the sign
along the path from $x$ to $y$. Hence $K_{n,j}(x,y)=s_xs_y$ on $C$.
Each block is therefore positive semidefinite, and so is
$K_{n,j}|_{S\times S}$.

Now set
$$
A_S=\frac1{J_n}\sum_{j\in I_S}K_{n,j}|_{S\times S}.
$$
Obviously, $A_S$ is positive semidefinite. Since each entry of
$K_{n,j}$ has absolute value at most one, for $x,y\in S$ we have
$$
|K_n(x,y)-A_S(x,y)|
\le \frac1{J_n}\sum_{j\notin I_S}|K_{n,j}(x,y)|
\le \frac{\binom{M}{2}}{J_n},
$$
where the sum is over the indices in the average defining $K_n$.
\end{proof}

Now, we are ready to define the state using ultrafilter. A \emph{free ultrafilter} $\omega$ on the infinite index set $I$ of
$(H_n)_n$ is a family of subsets of $I$ closed under finite
intersections and passage to supersets. It contains every cofinite
subset of $I$ and exactly one of $E$ and $I\setminus E$ for each
$E\subseteq I$. For a bounded sequence
$(a_n)_{n\in I}$, its \emph{ultralimit} 
$$a=\lim_{n\to\omega}a_n$$
is defined by
$$
\{n\in I:|a_n-a|<\varepsilon\}\in\omega
\qquad\text{for every }\varepsilon>0.
$$
The existsence of this limit is guaranteed by the the compactness of a closed disk containing the sequence, 
and it is obvious that the limit is unique. 

\begin{lemma}\label{lem:positive-functional}
Fix a free ultrafilter $\omega$ on the indices. The formula
\begin{equation}\label{eq:state}
\Phi_n(A)=\frac{1}{|U_n|}\sum_{x,y\in U_n} A_{xy}K_n(x,y),
\qquad
\Phi(A)=\lim_{n\to\omega}\Phi_n(A),\qquad A\in\C_u[X],
\end{equation}
defines a functional $\Phi$ on $\C_u[X]$ that extends uniquely
to a state on $\Cum(X)$ satisfying $$\Phi(p)=1, \quad \Phi(L_\sigma)=0.$$
\end{lemma}
\begin{proof}
Fix $A$ with propagation at most $R_A>0$. Each row and column has at
most $M=N_X(R_A)$ nonzero entries. At most $M|U_n|$ terms can contribute
to $\Phi_n(A)$, and each has absolute value at most
$\norm{A}_{\mathrm{red}}$, since $|K_n(x,y)|\le1$.
Thus $\Phi_n(A)$ is bounded in $n$ for each $A$, and the fixed
ultralimit defines a linear functional $\Phi$. Since $K_n(x,x)=1$,
we have $\Phi_n(1)=1$ for all $n$, hence $\Phi(1)=1$.

Fix $n$ and denote
$$S_z=\{x\in U_n:A_{zx}\ne0\}.$$
The row bound gives
$|S_z|\le M$, while the column bound shows that only finitely many
rows $z\in X$ can meet $U_n$. Expanding the matrix product gives rise to
$$
\Phi_n(A^*A)=\frac{1}{|U_n|}
\sum_{z\in X}\sum_{x,y\in S_z}\overline{A_{zx}}A_{zy}K_n(x,y).
$$
For every intermediate row index $z\in X$, the set $S_z$ is
contained in the ambient ball $B_X(z,R_A)$. Its cardinality is
therefore bounded by $M$, so Lemma~\ref{lem:small-kernel} applies
even when the vertices of $S_z$ are far apart in the graph metric.
Replacing the kernel on $S_z$ by its positive semidefinite
approximation makes the
inner sum nonnegative. The original sum is real since $K_n$ is real
and symmetric. The absolute value of the change is at most
$\epsilon_n(\sum_{x\in S_z}|A_{zx}|)^2$, where
$\epsilon_n=\binom{M}{2}/J_n$. By Cauchy--Schwarz inequality, we have that
$$
\sum_{z\in X}\Big(\sum_{x\in S_z}|A_{zx}|\Big)^2
\le M\sum_{x\in U_n}\norm{A\delta_x}_2^2
\le M|U_n|\norm{A}_{\mathrm{red}}^2.
$$
Consequently, 
$$\Phi_n(A^*A)\ge -M\epsilon_n\norm{A}_{\mathrm{red}}^2.$$
For each fixed $A$, the number $M$ is fixed while $J_n\to\infty$,
so the error tends to zero. The ultralimit preserves inequalities
between bounded real sequences, hence
$\Phi(A^*A)\ge0$.

For the algebraic GNS construction, quotient by the null space
of the positive form $\langle a,b\rangle_\Phi=\Phi(a^*b)$.
If $g=M_f$ is a diagonal contraction or $g=V_t$ is a partial
translation operator, then $1-g^*g=b^*b$ for a bounded diagonal
operator $b$. We may take $b=M_{\sqrt{1-|f|^2}}$ in the first case
and $b=1_{X\setminus\operatorname{dom}t}$ in the second. It follows that
$$
\Phi(a^*a)-\Phi((ga)^*(ga))=\Phi((ba)^*(ba))\ge0.
$$
Thus these left multiplications preserve the null space and act
contractively on the quotient. By the decomposition in
Section~\ref{sec:partial-translations}, the left action of
$A=\sum_i M_{f_i}V_{t_i}$ is bounded by $\sum_i\norm{f_i}_\infty$.
Left multiplication by $A^*$ is its adjoint for the form, since
$\Phi(a^*Ab)=\Phi((A^*a)^*b)$. Completing gives a bounded
$*$-representation $\pi$. Since $[1]$ is a unit vector,
$$
|\Phi(A)|=|\langle[1],\pi(A)[1]\rangle|
\le \norm{\pi(A)}\le\norm{A}_{\max}.
$$
The representation extends to $\Cum(X)$ by the definition of the
maximal norm, and its vector state at $[1]$ is the extension of $\Phi$.

The values of $K_n$ on the diagonal and on the edges determine
the evaluations in \eqref{eq:state} on $D$ and $B_\sigma$. On the diagonal
$K_n(x,x)=1$, so $\Phi_n(D)$ is the average degree of $H_n$.
For each edge $(x,y)$, the two ordered entries of $B_\sigma$
contribute $\sigma_n((x,y))^2=1$ each. Hence
$$
\Phi_n(D)=\Phi_n(B_\sigma)=\frac{2|E(H_n)|}{|U_n|},
\qquad
\Phi_n(p)=1.
$$
Taking the ultralimit gives $\Phi(p)=1$ and
$\Phi(L_\sigma)=\Phi(D)-\Phi(B_\sigma)=0$.
\end{proof}

\section{Proof of the main theorem}\label{sec:main-proof}
In this section, we shall prove our main theorem. Before proving, we shall recall the following theorem of Marcus--Spielman--Srivastava \cite{MSS} to estimate the maximal norm of the operator of the signed adjacency matrices.
\begin{theorem}[{\cite[Theorems 3.2, 3.7, 4.4 and 5.2]{MSS}}]\label{thm:mss-signing}
Let $H$ be a nonempty finite simple graph with degree at most $d\geq 2$. Then, there is a signing 
$$\sigma:E(H)\to\{-1,1\},$$
such that the largest eigenvalue of its signed adjacency matrix
$B_\sigma$ satisfies
$$
\lambda_{\max}(B_\sigma)\le 2\sqrt{d-1}.
$$
\end{theorem}

For the graphs constructed in the Lemma~\ref{lem:extraction}, we take $d=6$.
The above eigenvalue bound gives an upper bound on $B_\sigma$
in the concrete operator order. The minimum degree of vertices is five, and
$5>2\sqrt5$, so subtracting $B_\sigma$ from the degree diagonal ensures
a strictly positive lower bound for the signed Laplacian on
$\ell^2(U)$. 

\begin{theorem}\label{thm:converse}
Let $X$ be a discrete metric space of bounded geometry. If $X$ does not
have property A, then the quotient map
$\lambda:\Cum(X)\to\Cu(X)$ is not an isomorphism.
\end{theorem}
\begin{proof}
Suppose that $X$ does not have property A. Let $(H_n)_n$ be the sequence of graphs obtained in Lemmas~\ref{lem:witnesses}
and~\ref{lem:extraction}. By Theorem~\ref{thm:mss-signing}, for any a signing
$\sigma_n$ on each $H_n$, the signed adjacency matrix has largest
eigenvalue at most $2\sqrt5$. Let
$$
L_\sigma=D-B_\sigma,\qquad T=p-\frac{L_\sigma}{12}.
$$

Since the disjoint union of the $H_n$ has maximum degree six,
a greedy edge coloring uses at most eleven colors. Each color
class is a matching. Orient its edges and let $W_i$ be the resulting
partial translation operator weighted by the edge signs. Each
$W_i$ belongs to $\C_u[X]$ because its displacement is at most $L$.
Writing $q_i=W_iW_i^*$, we have $q_iW_i=W_i$. By summing over the
matchings, we have that
$$
D\pm B_\sigma
=\sum_i(q_i\pm W_i)^*(q_i\pm W_i).
$$
Therefore, $$0\le L_\sigma\le 2D$$ 
in $\Cum(X)$. The diagonal $6p-D$
has a square root in $\C_u[X]$, so $2D\le12p$, and hence
$$0\le T\le p.$$

In the reduced completion, the signing bound gives
$B_\sigma\le 2\sqrt5\,p$, while $D\ge5p$. It follows that
$$
0\le\lambda(T)\le c\,\lambda(p),
\qquad c=1-\frac{5-2\sqrt5}{12}<1.
$$
On the maximal completion, the state $\Phi$ in Lemma~\ref{lem:positive-functional} staisfies that 
$$\Phi(p)=1,~ \mbox{and}~\Phi(L_\sigma)=0.$$
Thus, $\Phi(T)=1$. It folllows from $0\le T\le p$ that 
$$\norm{T}_{\max}=1.$$

Set $a=(T-c1)_+$ by continuous functional calculus. Since $T$ is
positive and has maximal norm one, its spectrum in $\Cum(X)$
contains $1$, so $\norm{a}_{\max}=1-c>0$. On the other hand,
$0\le\lambda(T)\le c1$, and therefore
$$
\lambda(a)=(\lambda(T)-c1)_+=0.
$$
Thus $a$ is a nonzero element of $\ker\lambda$.
\end{proof}

Together with \cite[Proposition 1.3]{SW}, we have the following:
\begin{theorem}
Let $X$ be a discrete metric space with bounded geometry. Then, $X$ has proeprty A if and only if the canonical quotient map 
$$
\lambda: C^*_{u,max}(X)\to C^*_u(X)
$$
is an isomorphism. 	
\end{theorem}

\vspace{1cm}

\noindent Jintao Deng\\
\noindent Department of mathematics, Center of Noncommutative Geometry, Shanghai University of Finance and Economics\\
\noindent Email: \url{jintaodeng@sufe.edu.cn}
  
\vspace{1cm}

\noindent Geng Tian\\
\noindent Department of Mathematics, Liaoning University, Liaoning, China 110036\\
\noindent Email: \url{gengtian.ncg@gmail.com}

\end{document}